\documentclass{amsart}

\usepackage{amsmath,amsthm,amssymb,amscd,enumerate,mathtools,enumitem,tikz-cd,bm}
\usepackage{amsfonts}
\usepackage{rotating}
\usepackage{euscript}
\usepackage{pst-node}
\usepackage{epsfig,verbatim}

\def\be{\begin{equation}}
\def\ee{\end{equation}}

\def\C{{\mathbb C}}

\def\P{{\mathbb P}}

\def\R{{\mathbb R}} 
\def\Q{{\mathbb Q}}

\def\phi{{\varphi}}

\def\deg{{\rm deg\,}}

\def\bp{\begin{proposition}}
\def\ep{\end{proposition}}

\def\bt{\begin{theorem}}
\def\et{\end{theorem}}
\def\br{\begin{remark}}
\def\er{\end{remark}}
\def\be{\begin{equation}}
\def\bee{\begin{equation*}}
\def\l{\label}

\def\ee{\end{equation}}
\def\eee{\end{equation*}}
\def\bl{\begin{lemma}}
\def\el{\end{lemma}}
\def\bc{\begin{corollary}}
\def\ec{\end{corollary}}
\def\pr{\noindent{\it Proof. }}

\def\bd{\begin{definition}}
\def\ed{\end{definition}}

\def\h{\widehat}
\def\hat{\widehat}

\def\h{{\mathfrak{h}}} 

\newtheorem{theorem}{Theorem}[section]
\newtheorem{lemma}[theorem]{Lemma}
\newtheorem{definition}[theorem]{Definition}
\newtheorem{corollary}[theorem]{Corollary}
\newtheorem{proposition}[theorem]{Proposition}
\newtheorem{problem}[theorem]{Problem}

\theoremstyle{definition}

\theoremstyle{definition}
\newtheorem{remark}[theorem]{Remark}
\def\bpr{\begin{problem}}
\def\epr{\end{problem}}

\mathtoolsset{showonlyrefs}

\begin{document}

\title[]{Rational functions sharing preimages  and \linebreak height functions
}

\author[F. Pakovich]{Fedor Pakovich}
\thanks{
This research was supported by ISF Grant No. 1092/22} 
\address{Department of Mathematics, 
Ben Gurion University of the Negev, P.O.B. 653, Beer Sheva,  8410501, Israel}
\email{pakovich@math.bgu.ac.il}

\begin{abstract} 
Let \( A \) and \( B \) be non-constant rational functions over $\C$, and let \( K \subset \mathbb{P}^1(\mathbb{C}) \) be an infinite set. Using height functions, we prove that the inclusion  \( A^{-1}(K) \subseteq B^{-1}(K) \) implies the inequality \( \deg B \geq \deg A \) in the following two cases: the set \( K \) is contained in \( \mathbb{P}^1(\bm{k}) \), where \( \bm{k} \) is a finitely generated subfield of \( \mathbb{C} \), or the set \( K \) is discrete in \( \mathbb{C} \), and \( A \) and \( B \) are polynomials. In particular, this implies that for $A$, $B$, and \( K \) as above, the equality \( A^{-1}(K) = B^{-1}(K) \) is impossible, unless \( \deg A = \deg B \).
 
 \end{abstract}


\maketitle

\section{Introduction} 
Let \( A \) and \( B \) be polynomials over \( \mathbb{C} \), and let \( K \subset \mathbb{C} \) be a compact set, either finite or infinite. The problem of characterizing collections \( A \), \( B \), and \( K \) such that \( A \) and \( B \) share the preimages of \( K \), i.e.,  
\be \l{00} 
A^{-1}(K) = B^{-1}(K), 
\ee  
was solved in a series of papers \cite{d}, \cite{d2}, \cite{p1}, \cite{p2}, \cite{p3}. Notice that this problem can be viewed as describing compact, completely invariant sets of polynomial correspondences. It also includes, as a special case, the classification of polynomials that share Julia sets (see \cite{a1}, \cite{b}, \cite{be1}, \cite{be2}, \cite{f}, \cite{sh}).

The results of  \cite{d}, \cite{d2}, \cite{p1}, \cite{p2}, \cite{p3} can be summarized as follows. It was shown in \cite{p1}, \cite{p2} that for any compact set \( K \subset \mathbb{C} \) containing at least two points and for polynomials \( A \) and \( B \) of the same degree, the equality \eqref{00} implies that \( A = \sigma \circ B \) for some polynomial \( \sigma \) of degree one satisfying \( \sigma(K) = K \).  
For polynomials of arbitrary degrees, a description of solutions to \eqref{00} for infinite compact sets \( K \subset \mathbb{C} \) 
was given 
in \cite{d}, \cite{d2}. Namely, it was shown that if \( K \) is not a union of circles or a segment and \( \deg B \geq \deg A \), then there exists a polynomial \( G(z) \) such that \( B = G \circ A \) and \( G^{-1}(K) = K \).

A description of solutions to a more general set-theoretic equation  
\be \l{0} A^{-1}(K_1) = B^{-1}(K_2),\ee  
where \( A \) and \( B \) are polynomials and \( K_1 \) and \( K_2 \) are arbitrary compact subsets of \( \mathbb{C} \), not necessarily equal, was obtained in \cite{p3}. Specifically, in \cite{p3}, condition \eqref{0} was related to the functional equation  
\begin{equation}  
\label{2}  
G \circ A = F \circ B,  
\end{equation}  
in polynomials. It is clear that for any polynomial solution of \eqref{2} and any compact set \( K \subset \mathbb{C} \), one obtains a solution of \eqref{0} by setting  
\begin{equation}  
\label{3}  
K_1 = G^{-1}(K), \quad K_2 = F^{-1}(K),  
\end{equation}  
and the main result of \cite{p3} states that {\it all} solutions to \eqref{0} can be constructed in this way, provided that the compact set defined by either side of \eqref{0} contains at least \( \text{LCM}(\deg A, \deg B) \) points. 
Since solutions to \eqref{2} are characterized by Ritt’s theory of polynomial factorization \cite{r1}, this provides a rather precise description of solutions to \eqref{0}, which can then be applied to various related problems  (see \cite{p3} and \cite{pj} for further details).

The methods employed in the aforementioned papers are restricted to the polynomial case, and the problem of finding solutions to \eqref{00} and \eqref{0} when \( A \) and \( B \) are \textit{rational} functions and \( K_1 \) and \( K_2 \) are compact subsets of \( \mathbb{C} \) or $\P^1(\C)$ remains largely unresolved. It appears that the only significant result in this direction is due to Bellaïche (\cite{bel}), who investigates the broader problem of describing invariant sets of correspondences. In this note, we consider a particular class of correspondences defined by a pair of rational functions \( A \) and \( B \) on \( \mathbb{P}^1(\mathbb{C}) \). For such a correspondence \( (\mathbb{P}^1, A, B) \), its forward and backward maps are given by the multivalued functions \( B(A^{-1}(z)) \) and \( A(B^{-1}(z)) \), respectively. In this context, the problem of describing the sets \( K \) satisfying \eqref{00} is clearly equivalent to characterizing the completely invariant sets of the correspondence \( (\mathbb{P}^1, A, B) \), that is, the sets stable under both the forward and backward maps of \( (\mathbb{P}^1, A, B) \).

When applied to equation \eqref{00}, the main result of \cite{bel} implies that if for rational functions $A$ and $B$ equality \eqref{00} holds for infinitely many finite sets \( K \), then there exists a rational function \( F \) such that  
\be \label{4}  
F \circ A = F \circ B.  
\ee  
Notice that in the last case \eqref{00} is satisfied for every set of the form \( K = F^{-1}(\hat{K}) \), where \( \hat{K} \subset \P^1(\C) \). Notice also that equality \eqref{4}  obviously implies the equality $\deg A=\deg B$.

In this note, we consider solutions to \eqref{00}, where \( A \) and \( B \) are arbitrary non-constant rational functions. However, instead of assuming that \( K \) is compact, we assume that \( K \) is an infinite set contained in \( \mathbb{P}^1(\bm{k}) \), where \( \bm{k} \) is a finitely generated subfield of \( \mathbb{C} \)—for instance, a number field. 
In fact, instead of condition \eqref{00}, we consider the more general condition  
\be \l{1}  
A^{-1}(K) \subseteq B^{-1}(K),  
\ee  
and establish the following result.

\begin{theorem} \l{t1}  Let \( A \) and \( B \) be non-constant rational functions over $\C$, and let \( K\) be an infinite subset of \( \mathbb{P}^1({\bm k}) \), where \( \bm{k} \) is a finitely generated subfield of \( \mathbb{C} \). If \( A^{-1}(K) \subseteq B^{-1}(K) \), then \( \deg B \geq \deg A \). In particular, $A^{-1}(K) \neq  B^{-1}(K)$, unless $\deg A = \deg B$. 
\end{theorem}

As an immediate corollary, we obtain the following.

\bc \l{c1}  
Let $A$ and $B$ be non-constant rational functions over $\C$. Then the correspondence $(\C\P^1, A,B)$ can have an infinite completely invariant set  contained in  a finitely generated subfield of $\mathbb{C}$ only if $\deg A = \deg B$. 
\ec

We also prove a result concerning  solutions to  \eqref{1}, where \( A \) and \( B \) are polynomials, and \( K \) satisfies a condition that is, in a sense, opposite to compactness---discreteness.  
We show that the same conclusion as in Theorem \ref{t1} holds in this setting as well.

\begin{theorem} \l{t2}  
Let \( A \) and \( B \) be non-constant polynomials over $\C$, and let \( K \) be an infinite discrete subset of \( \mathbb{C} \). If \( A^{-1}(K) \subseteq B^{-1}(K) \), then \( \deg B \geq \deg A \).   In particular, $A^{-1}(K) \neq  B^{-1}(K)$, unless $\deg A = \deg B$. 
\end{theorem}  

Again, the following corollary is immediate.

\begin{corollary} \l{cc2}  
Let \( A \) and \( B \) be non-constant polynomials over $\C$. Then the correspondence $(\C\P^1, A,B)$ can have an infinite completely invariant set that is discrete in $\C$ only if  $\deg A = \deg B$. 
\end{corollary}

Notice that examples of  \( A \) and \( B \) of equal degree satisfying \eqref{00} for sets \( K \) as described above do exist. For example, if \( K \) is the set of squares of integers, then for the polynomials  
\[
A = z^2, \quad B = (z+1)^2,   
\]
we clearly have  
\[
A^{-1}(K) = B^{-1}(K) = \mathbb{Z}.
\]
Here, \( K \) is both discrete in \( \mathbb{C} \) and contained in $\Q$. Let us mention that, in contrast to the case where \( K \) is compact, the polynomial \( A \) cannot be represented as \( \sigma \circ B \) for some polynomial \( \sigma \) of degree one, since \( A \) and \( B \) have different critical points. The equality \eqref{4} is also impossible because, if it held, the rational function defined by any part of this equality would be invariant under \( z \to z+1 \).

Our proof of Theorem \ref{t1} is similar to the proof, based on the Weil height, of the well-known fact that the inverse orbit of a point \( x\in \overline {\Q} \) under a rational function \( F\in\overline{\Q}(z) \) can contain only finitely many points in any finite extension of \( \Q \). However, instead of the orbit of a rational function, we consider the orbit of the algebraic function \( B(A^{-1}(z)) \), and instead of the Weil height on $\P^1(\overline {\Q})$ we use the Moriwaki height  on $\P(\overline{\bm k}).$ The use of the Moriwaki height eliminates any restrictions on the coefficients of frunctions \( A \) and \( B \), allowing Theorem \ref{t1} and Corollary \ref{c1} to apply to a wide variety of infinite sets \( K \). 
For example, \( K \) can be a subset of a number field, a subset of the field \( \mathbb{Q}(\pi, e) \), a subset of an orbit of a rational function, etc.

Notice that if \( K \) is a compact subset of \( \mathbb{C} \mathbb{P}^1 \), the condition \eqref{1} does not generally imply the inequality \( \deg B \geq \deg A \), since there exist rational functions of different degrees satisfying \eqref{00}. Indeed, it is enough to take any rational function \( P \) and a set \( K \) satisfying \( P^{-1}(K) = K \), and set \( A = Q \), \( B = P \circ Q \), where \( Q \) is any rational function.

The proof of Theorem \ref{t2} is somewhat similar to that of Theorem \ref{t1}, as we use the function \(\h(z)= \log(\max\{1, \vert z \vert\}) \) as an analogue of the height function. However, instead of seeking a contradiction with the Northcott property, which is not present in this setting, we instead deduce the contradiction from the absence of limit points. In both cases, however, the key point is the inequality  
$$
\big| \mathfrak{h}(R(z)) - \deg R \cdot \mathfrak{h}(z) \big| < C,  
$$  
where \( C = C(R) \) is a constant, which holds for any \( R \) in a suitable subset of \( \C(z) \) and any \( z \) in a suitable subset of \( \P^1(\C) \).

\section{Proofs}

For non-constant rational functions $A$ and $B$ over $\C$, we define the multi-valued function ${H}={H}_{A,B}$ by the formula $${H}(z) = B(A^{-1}(z)).$$ More precisely, for a point $z \in \P^1(\C)$ or, more generally, for a subset $K$ of $\P^1(\C)$, the notation ${H}(z)$ or ${H}(K)$ refers to the set $B(A^{-1}(z))$ or $B(A^{-1}(K))$, respectively. 
Inductively, we define the set ${H}^{\circ k}(K)$ for $k \geq 2$ by  
\[
{H}^{\circ k}(K) = {H}({H}^{\circ (k-1)}(K)), \quad k \geq 2.
\]  
Finally, we call the set  
\[
\bigcup_{i=1}^{\infty} {H}^{\circ i}(K)
\]  
the orbit of $K$ under ${H}$.

To establish our results, we make use of the existence
of a finite set \( K \) whose orbit under \( H \) is infinite.  
Notice that such a set \( K \) may not always exist. For instance, if \( A \) and \( B \) satisfy \eqref{4} for some rational function \( F \), then it is easy to see that for any set \( K \), the orbit of \( K \) under \( H \) is contained in \( F^{-1}(F(K)) \).  
Nevertheless, if \( \deg A \neq \deg B \), such sets always exist, as the following elementary lemma shows.

\bl \l{len} Let $A$ and $B$ be rational functions over $\C$ of degrees $n$ and $m$ with $n > m>0$, and let $K$ be a finite subset of $\mathbb{P}^1(\mathbb{C})$. Then the orbit of $K$ under ${H}$ is infinite whenever  
\be\l{in}  
\vert K\vert > \frac{2n-2}{n-m}  
\ee  
\el  
\pr It  follows easily from the Riemann-Hurwitz formula that for a rational function $F$ of degree $d$ and a finite subset $K$ of $\mathbb{P}^1(\mathbb{C})$, the inequality  
\[
 \vert F^{-1}(K)\vert  \geq d(\vert K\vert-2)+2
\]  
holds. On the other hand, setting $K = F(T)$ in an obvious inequality   
\[
\vert F^{-1}(K)\vert \leq d\vert K\vert,  
\]  
 we obtain that for any finite subset $T$ of $\mathbb{P}^1(\mathbb{C})$, the inequality  
\[
\vert F(T) \vert \geq \frac{\vert F^{-1}(F(T)) \vert}{d} \geq \frac{\vert T \vert}{d}  
\]  
holds.  

These inequalities imply that  
\be \l{if}  
\vert {H}(K)\vert \geq \frac{\vert A^{-1}(K)\vert}{m} \geq \frac{n(\vert K\vert-2)+2}{m}.  
\ee  
Since inequality \eqref{in}   
implies   the inequality  
\[
n(\vert K\vert-2)+2> m\vert K\vert,  
\]  
it follows from \eqref{if} that $ \vert {H}(K)\vert >\vert K\vert $ whenever \eqref{in} holds. Applying this inequality recursively, we conclude that the orbit of $K$ under ${H}$ is infinite.  
\qed

Let $A$ and $B$ be non-constant rational functions over $\C$, and let \( k \) be an infinite subset of \( \mathbb{P}^1(\mathbb{C}) \)  such that $  A^{-1}(k)\subseteq k .$  
Assume that there exist a function \linebreak $\h: k\rightarrow \R_{\geq 0}$ and constants $C_1,C_2> 0$ such that 
for every $z\in k$ the inequalities  
\be \l{proo} \vert \mathfrak{h}(A(z))-\deg A\cdot\mathfrak{h}(z)\vert < C_1, \ \ \ \vert \mathfrak{h}(B(z))-\deg B\cdot \mathfrak{h}(z)\vert < C_2\ee hold.

\bp \l{p1} Under the above assumptions on \( k \), $A$, $B$, and \( \mathfrak{h} \), suppose that there exists an infinite subset \( K \subseteq k \) such that \( A^{-1}(K) \subseteq B^{-1}(K) \) and \linebreak\( \deg A > \deg B \). Then there exists a constant \( M > 0 \) such that the inequality  
$\h(z) < M  $
holds for infinitely many \( z \in K \).
\ep 
\pr We apply inequalities \eqref{proo} to the orbit \( H(z) \) with \( z \in K \), noting that the condition \( A^{-1}(K) \subseteq B^{-1}(K) \) ensures that such an orbit remains in \( K \). Let us set \( n = \deg A \) and \( m = \deg B \). For \( z \in K \), substituting any \( z' \in A^{-1}(z) \) into the first inequality in \eqref{proo}, we obtain  
\[
\frac{\mathfrak{h}(z)}{n} - \frac{C_1}{n} < \mathfrak{h}(z') < \frac{\mathfrak{h}(z)}{n} + \frac{C_1}{n}.
\]  
Next, applying the second inequality in \eqref{proo}, we find that for any \( \hat{z} \in H(z) \), the inequality  
\[
m \mathfrak{h}(z') - C_2 < \mathfrak{h}(\hat{z}) < m \mathfrak{h}(z') + C_2
\]  
holds. Thus,  
\begin{equation} \l{bew} 
\left(\frac{m}{n}\right) \mathfrak{h}(z) - C < \mathfrak{h}(\hat{z}) < \left(\frac{m}{n}\right) \mathfrak{h}(z) + C, 
\end{equation}
where \( C = \frac{m}{n} C_1 + C_2 \).

Applying  \eqref{bew} recursively, we conclude that if $z$ is a point of $K$, then  
for any point $\hat z\in {H}^{\circ k}(z)$ the inequality 
\be \l{fu} \mathfrak{h}(\hat z)<
\left(\frac{m}{n}\right)^k\mathfrak{h}(z)+C\left(\left(\frac{m}{n}\right)^{k-1}+\left(\frac{m}{n}\right)^{k-2}+\dots +1
\right)<\mathfrak{h}(z)+\frac{C}{1-\frac{m}{n}}\ee
holds. 
Furthermore, if $K$ is a finite subset of $k$ such that inequality \eqref{in} holds and 
$M=\max_{z\in K}\mathfrak{h}(z),$ then \eqref{fu} implies that for every point $\hat z\in {H}^{\circ k}(K)$ we have 
$$\mathfrak{h}(\hat z)<M+\frac{C}{1-\frac{m}{n}}.$$ Since the orbit of $K$ is infinite by Lemma \ref{len}, this completes the proof. \qed  

Theorem \ref{t1} and \ref{t2} are obtained from Proposition \ref{p1} by using appropriate $k$ and  $\h$.
To prove Theorem \ref{t2}, we set $k=\C$ and define $\h$ by the formula  
$$\h(z)=\log(\max\{1, \vert z \vert\}).$$ It is well known that for every non-constant polynomial $R$ there exists $C>0$ such that  the inequality 
\be \l{pr}\big| \mathfrak{h}(R(z)) - \deg R\cdot \mathfrak{h}(z) \big| < C \ee 
holds for every $z\in \C$. 
Thus,  since $R^{-1}(\C)=\C$, Proposition \ref{p1} is applicable.

\vskip 0.2cm

\noindent {\it Proof of Theorem \ref{t2}.} It follows from  Proposition \ref{p1} that if $\deg A>\deg B$, then there exists $M>0$ such that the inequality $\vert z\vert <M$ holds for infinitely many $z$ in $K$. Since $K$ is infinite, this implies that $K$ has a limit point in $\C$. \qed 

\vskip 0.2cm 

To prove Theorem \ref{t1}, we use the Moriwaki height, which generalizes the Weil height. Let us recall that one can define the (logarithmic) Weil height $h (x)$ on points of \( \mathbb{P}^1(\overline{\Q}) \), satisfying the following two properties (see, e.g., \cite{silver}). First, for any rational function \( R \in \overline{\Q}(z) \), there exists $C > 0$ such that \eqref{pr} holds for every $z\in \overline{\Q}$. Second, for any numbers \( D_1, D_2 > 0 \), there are only finitely many points \( x \in \mathbb{P}^1(\overline{\Q}) \) satisfying the conditions  
$$ h(x) \leq D_1, \quad [\mathbb{Q}(x):\mathbb{Q}] \leq D_2 $$  
(Northcott's theorem).  
Furthermore, it follows from the results of Moriwaki (\cite{mor}) that one can define a height function \( \h \) with similar properties on \( \mathbb{P}^1(\overline{\bm{k}}) \)  for any field \( \bm{k} \) finitely generated over \( \mathbb{Q} \). 

Referring the reader to \cite{tits} for more detail on the Moriwaki height in the considered context, we state only properties of $\h$ needed for our purpose (see \cite{mor}, \cite{tits}). First, if \( \h \) is the Moriwaki height on \( \mathbb{P}^1(\overline{\bm{k}}) \), then for any rational function \( R \in {\bm{k}}(z) \), there exists \( C > 0 \) such that \eqref{pr} holds for every \( z \in \overline{\bm{k}} \). Second, for any \( D_1, D_2 > 0 \), only finitely many points \( x \in \mathbb{P}^1(\overline{\bm{k}}) \) satisfy the conditions  
\[
\h(x) \leq D_1, \quad [\bm{k}(x):\bm{k}] \leq D_2.
\]

\vskip 0.2cm

\noindent {\it Proof of Theorem \ref{t1}.} Adjoining to $\bm k$ coefficients of $A$ and $B$, without loss of generality we may assume that $A,B\in {\bm k}(z)$, and it clear that $A^{-1}(\overline{{\bm k}})\subset \overline{{\bm k}}$. Setting now $k=\overline {\bm k}$ and using  Proposition \ref{p1} we see that if $\deg A>\deg B$, then there exists $M>0$ such that the inequality $\h(z) <M$ holds for infinitely many $z$ in $K$ in contradiction with the Northcott property, as $K\subseteq \bm k.$  \qed

\end{document}